\documentclass[12pt]{amsart}

\usepackage[centertags]{amsmath}
\usepackage{amsfonts}
\usepackage{amscd}
\usepackage{latexsym}
\usepackage{amssymb}
\usepackage{mathrsfs}
\usepackage{euscript}
\usepackage[all]{xy}
\usepackage{pxfonts}

\newcommand{\B}[1]{{\mathbb #1}}

\newtheorem{theorem}[subsection]{Theorem}

\newtheorem{corollary}[subsection]{Corollary}
\newtheorem{lemma}[subsection]{Lemma}
\newtheorem{proposition}[subsection]{Proposition}
\theoremstyle{definition}
\newtheorem{definition}[subsection]{Definition}
\newtheorem{example}[subsection]{Example}
\theoremstyle{remark}
\newtheorem{remark}[subsection]{Remark}
\newtheorem{question}[subsection]{Question}

\numberwithin{figure}{section}
\numberwithin{table}{section}
\numberwithin{equation}{section}

\newcommand{\al}{{\alpha}}

\newcommand{\be}{{\beta}}

\newcommand{\Om}{{\Omega}}
\newcommand{\om}{{\omega}}

\newcommand{\De}{{\Delta}}
\newcommand{\ga}{{\gamma}}
\newcommand{\Ga}{{\Gamma}}

\newcommand{\si}{{\sigma}}
\newcommand{\Si}{{\Sigma}}

\newcommand{\Mo}{(M,\omega )}
\newcommand{\Wo}{(W,\omega_W )}

\newcommand\SL{\operatorname{SL}}

\newcommand\Diff{\operatorname{Diff}}
\newcommand\Flux{\operatorname{Flux}}
\newcommand\Symp{\operatorname{Symp}}
\newcommand\Ham{\operatorname{Ham}}

\newcommand\vol{\operatorname{vol}}

\newcommand{\ms}{{\medskip}}
\newcommand{\bs}{{\bigskip}}
\newcommand{\NI}{{\noindent}}

\begin{document}

\title{Symplectically hyperbolic manifolds}
\author{Jarek K\c edra}
\address{Mathematical Sciences, University of Aberdeen, 
         Meston Building, Aberdeen, AB243UE, Scotland, UK\\
Institute of Mathematics, University of Szczecin,
Wielkopolska 15, 70-451 Szczecin, Poland}
\email{kedra@maths.abdn.ac.uk}
\urladdr{http://www.maths.abdn.ac.uk/\~{}kedra}
\date{\today }
\thanks{thanks} 
\subjclass[2000]{Primary 53; Secondary 57}
\begin{abstract}
A symplectic form is called hyperbolic if its
pull-back to the universal cover is a differential of 
a bounded one-form. 
The present paper is concerned with the properties and
constructions of manifolds admitting hyperbolic symplectic
forms.
The main results are:
\begin{itemize}
\item
If a symplectic form represents a bounded
cohomology class then it is hyperbolic. 
\item
The symplectic hyperbolicity is equivalent
to a certain isoperimetric inequality. 
\item
The fundamental group of symplectically hyperbolic manifold
is non-amenable.
\end{itemize}
We also construct hyperbolic symplectic forms on certain
bundles and Lefschetz fibrations, discuss the dependence
of the symplectic hyperbolicity on the fundamental group and
discuss some properties of the group of symplectic diffeomorphisms
of a symplectically hyperbolic manifold.

\setcounter{tocdepth}{1}
\tableofcontents

\end{abstract}


\maketitle

\section{Introduction and statements of the results}\label{S:intro}

Let $\om \in \Om^k(M)$ be a closed differential form on
a closed manifold~$M$. Let $p:\widetilde M\to M$ be
the universal covering. Let $g$ be a Riemannian metric
on $M$ and $\tilde g$ the induced metric on the
universal cover. 

\begin{definition}[Gromov \cite{MR1085144}]\label{D:d-bounded}
The form $\om$ is called
{\em $\tilde d$-bounded} if its pull-back is a
differential of a bounded form. That is  
$p^*\om=d\alpha$ and there exist a constant
$C\in \B R$ such that $\sup_{x\in M}|\alpha(x)|\leq C$.
The norm of a differential form is defined by
~$|\om(x)|:= \max \{\om(x)(X_1,\ldots,X_k)\,|\,\,|X_i|=~1\}$.

\end{definition}

\begin{proposition}\label{P:cohomology_class}
The $\tilde d$-boundedness does not depend on the
choice of a Riemannian metric on $M$. Moreover, if $\om $ is
$\tilde d$-bounded then so is $\om + d\xi$. 
\end{proposition}

\begin{proof}
Since $M$ is compact, the first statement is clear.
If $p^*\om =d\alpha $ then $p^*(\om +d\xi) = d(\alpha + p^*(\xi))$.
Again, due to the compactness of $M$, $\alpha $ is bounded 
if and only if $\al + p^*(\xi)$ is bounded.
\end{proof}

\begin{definition}\label{D:hyperbolic}
Let $\Mo$ be a closed symplectic manifold. If the symplectic
form is $\tilde d$-bounded then it is called {\bf hyperbolic}.
A manifold $\Mo$ is then called {\bf symplectically hyperbolic}.
\end{definition}

\subsection*{Aim}
The present paper provides some constructions of
symplectically hyperbolic manifold and investigates 
their geometric and topological properties.

\subsection*{History}\label{SS:history}
The $\tilde d$-boundedness was defined by Gromov in
\cite{MR1085144}, where he proved a version of the 
Lefchetz theorem for $L^2$-cohomology for a symplectically
hyperbolic manifold $\Mo$ which is K\"ahler.
As a corollary he obtained that $(-1)^n\chi(M)> 0$,
where $\dim M=2n$. This is a particular case of a
conjecture (attributed to Hopf) stating that if an even 
dimensional manifold $M$ is negatively curved then its 
Euler characteristic satisfies $(-1)^n\chi(M)> 0$.

Properties of $\tilde d$-bounded forms were also
investigated by Sikorav ~\cite{MR1871291}, where he proved
certain isoperimetric inequalities.

Polterovich in \cite{MR2003i:53126} proved a 
number of results about symplectic diffeomorphisms of 
symplectically hyperbolic manifolds. 
For example, he proved that there are
strong restrictions on finitely generated groups 
which admit a Hamiltonian representation on a symplectically
hyperbolic manifolds. For example $\SL(n,\B Z)$ does not
admit any such representation. A Hamiltonian representation
of a group $G$ o a symplectic manifold $\Mo$
is a homomorphism $G\to \Ham\Mo$.

\subsection*{Basic examples and properties}
\begin{example}\label{E:surface}
Let $(\Si,\om)$ be a closed surface of genus at least 2.
Let $g$ be a hyperbolic metric (i.e. the sectional curvature equal
to $-1$). The universal cover is the hyperbolic
plane $\B H=\{(x,y)\in \B R^2\,|\, y>0\,\}$. Let 
$\om $ be chosen so that the induced form on the universal
cover is equal to $\frac{1}{y^2}\cdot dy\wedge dx$.
We have that $\widetilde \om= d(\frac{dx}{y})$.
The hyperbolic metric on $\B H$ is given by $\frac 1 {y^2} (dx^2+dy^2)$ 
so we calculate 
$$\left |\frac{dx}{y}\right|=
\max \left \{\frac{dx}{y}(a\partial _x + b\partial _y)\,\big|\, 
a^2+b^2=y^2 \right \}=
\max \left\{\frac{a}{y}\,\big|\, a^2+b^2=y^2\right\}= 1$$
for any point $(x,y)\in \B H$.
\end{example}

\begin{example}\label{E:product}
The product of symplectically hyperbolic manifolds
is symplectically hyperbolic. Hence, it follows from
the previous example that the product form on a
product of surfaces 
$\Si_1\times \dots \times \Si_k$ of genus at least 2 
is hyperbolic. 

Since such a product contains a lagrangian torus the product
symplectic form can be slightly perturbed to a form
which does not vanish on this torus. Thus it is
not hyperbolic anymore, due to Proposition
\ref{P:properties} (cf. Example \ref{E:lattice}).
\end{example}

\begin{example}\label{E:immersion}
Let $\Mo$ be a symplectically hyperbolic and let
$f:S\to M$ be a symplectic immersion. Then
$f^*\om$ is hyperbolic. We have the diagram of
the universal covers.
$$
\xymatrix
{
\widetilde S \ar[d]^{p_S}\ar[r]^{\tilde f} & \widetilde M\ar[d]^{p_M}\\
S \ar[r]^f & M
}
$$ 
Let $p_M^*\om = d\alpha$, where $\alpha $ is a bounded one-form.
Let $S$ be equipped with the Riemannian metric induced form
$M$. Then we have that $p_S^*f^*\om=\tilde f^*p_M^*\om =
\tilde f^*d\alpha = d(\tilde f^*\alpha)$ and 
$\tilde f^*\alpha $ is bounded with respect to the
above mentioned metric.
\end{example}

\begin{example}
Let $f:\Mo\to \Wo$ be a smooth map between symplectic
manifolds such that $f^*[\om_W]=[\om]$. If $\om_W$ is
hyperbolic then so is $\om$, according to Proposition
\ref{P:cohomology_class}. In particular, branched covers
of symplectically hyperbolic manifolds are symplectically
hyperbolic (see Gompf \cite{MR99k:57068} for more detailed 
examples and calculations).

\end{example}

\begin{example}\label{E:torus}
The torus $T^n$ equipped with the standard symplectic
form is not symplectically hyperbolic. Let $g$ be the
standard flat metric. The universal cover is the
standard flat $\B R^{2n}$. Let $B(r)\subset \B R^{2n}$ denote
the ball of radius $r$. For any primitive $\alpha $ 
we calculate
\begin{eqnarray*}
\frac{(2\pi)^n}{2\cdot4\cdots (2n)}\cdot r^{2n} 
&=& \vol(B(r))=\int_{B(r)}\widetilde \om 
=\int_{\partial B(r)} \alpha \\
&\leq& 
\vol(\partial B(r))\cdot \sup_{x\in B(r)}|\alpha(x)|\\
&=&
\frac{(2\pi)^n}{2\cdot 4\cdots (2n-2)}\cdot
r^{2n-1}\cdot \sup_{x\in B(r)}|\alpha(x)|
\end{eqnarray*}
Hence we get that
$$\sup_{x\in B(r)}\,|\alpha(x)| \geq \frac r{2n}$$
which proves that every primitive of $\widetilde \om$
is unbounded.
\end{example}

\begin{proposition}\label{P:properties}
If $\om$ is a hyperbolic symplectic form then
$$
\int_{f(\Si)}\om = 0
$$
for any map $f:\Si\to M$, where $\Si$ is either a sphere
or a torus. In particular, the same statement holds if
$[\om]$ is bounded.
\end{proposition}

\begin{proof}
If $\om$ did not vanish on a sphere $f:S^2\to M$ then
the induced form $\widetilde \om$ would be non-zero in
the cohomology of the universal cover. This contradicts
the definition.

Let $f:T^2\to M$ be any smooth map. Since the pull-back 
$f^*[\om]\neq 0$, the standard symplectic form is
hyperbolic which contradicts Example \ref{E:torus}.
\end{proof}

\begin{question}
Suppose that a symplectic form vanishes on spheres and tori.
Is it hyperbolic?
\end{question}

\subsection*{A sufficient condition for symplectic hyperbolicity}
A cochain $\al\in C^*(X)$ on a topological space $X$ is
called {\bf bounded} if there exist a constant $C>0$ such that
$\left<a,s\right> < C$
for every singular simplex $s:\De\to X$.
A cohomology class is called {\bf bounded}  
if it is represented by a bounded
cochain (see Gromov \cite{MR686042}). 
The following lemma is a particular case of Theorem
\ref{T:bounded} proved in Section \ref{S:bounded}.

\begin{lemma}\label{L:bounded}
If a symplectic form represents a bounded cohomology class 
then it is hyperbolic.\qed
\end{lemma}

Many of the results of the present paper use this lemma.
That is, we construct symplectic forms which represent
bounded cohomology classes. It is known that a cohomology
class of degree at least two of a non-elementary hyperbolic group is
bounded. Hence we get the following easy application
of the above lemma.

\begin{corollary}\label{C:hyperbolic}
Let $\Mo$ be a closed symplectic manifold. If $[\om]$ is
aspherical and $\pi_1(M)$ is hyperbolic then $\om$
is hyperbolic. In particular, if $M$ admits a Riemannian
metric of negative sectional curvature then $\om $ is
hyperbolic. \qed 
\end{corollary}

\begin{remark}\hfill
~\begin{enumerate}
\item
A cohomology class $\al \in H^k(X)$
is called {\bf aspherical} if $\left <\al,f_*[S^k]\right> \\
=~0$
for any continous map $f:S^k\to X$. It is equivalent to
the fact that $\al=c_X^*(\Om)$, where 
$c_X:X\to K(\pi_1(X),1)$ is the classifying map
(see \cite{MR2103543,SAS} for more results about the
topology of manifolds admitting aspherical symplectic
forms).
\item
A hyperbolic group is called non-elementary if it 
does not contain a finite index cyclic group.
Since the fundamental group of a symplectic manifold
with an aspherical symplectic form is of virtual
cohomological dimension at least two, $\pi_1(M)$
in the above corollary is automatically non-elementary.
\item
There are no examples of closed manifolds of constant
negative sectional curvature admitting a symplectic form.
It is conjectured that such manifolds do not exist
\cite{MR1919897}.
\item
I do not know any example of a symplectic manifold
$\Mo$ such that $\om$ is hyperbolic and $[\om]$ is not
represented by a bounded cocycle.
\end{enumerate}
\end{remark}

\begin{example}\label{E:lattice}
Let $\B C\B H^n$ be a complex hyperbolic space. Its sectional
curvature is negative and pinched between $-1$ and $-1/4$.
The isometry group if isomorphic to $PU(n,1)$ \cite{MR1695450}. 
Thus any cocompact lattice in $PU(n,1)$ gives rise to a negatively
curved closed K\"ahler manifold. We get this way aspherical
symplectically hyperbolic manifolds whose fundamental groups
are hyperbolic and have trivial the first Betti number
\cite{MR1090825}.
Notice that any symplectic form is hyperbolic in this case
according to Corollary \ref{C:hyperbolic} 
(cf. Example \ref{E:product}).
\end{example}

\begin{question}
What is the relation between bounded and $\tilde d$-bounded
cohomology classes on a closed manifold?
\end{question}

\subsection*{Linear isoperimetric inequality}
We characterize symplectically hyperbolic manifolds by
a geometric condition inspired by a characterization 
of hyperbolic groups by certain isoperimetric inequality. 
More precisely, let $\Mo$ be a symplectic manifold and $g$ a 
Riemannian metric on $M$.
We say that the symplectic form $\om$ satisfies the
{\bf linear isoperimetric inequality} (with respect to $g$)
if there exists a constant
$C > 0$ such that for all $f:D^2\to M$ we have
$$
\int_{f(D^2)}\om \leq C\cdot \text{Length}(f(\partial D^2)).
$$

\begin{theorem}\label{MT:isoperimetric}
Let $\Mo$ be a closed symplectic manifold.
The symplectic form $\om $ satisfies a linear isoperimetric inequality
if and only if it is hyperbolic.
\end{theorem}

\subsection*{Fundamental group}
In Section \ref{S:group}, we prove that symplectic
hyperbolicity depends (in an appropriate sense)
on the fundamental group only.
More precisely, let $\Mo$ and $\Wo$ be closed symplectic
manifolds with isomorphic fundamental groups.
If $[\om]=c_M^*(\Om)$ and $[\om_W]=c_W^*(\Om)$ for a
class $\Om\in H^2(\pi_1(M);\B R)$ then $\om$ is
hyperbolic if and only if so is $\om_W$. Here,
$c_X:X\to K(\pi_1(X),1)$ denotes the map classifying
the universal cover.

Recall that the space $K(G,1)$ is usualy not a manifold so
the $\tilde d$-boundedness cannot be directly generalized.
The above result shows that the notion of symplectic hyperbolicity
makes sense for classes $\Om\in H^2(G,\B R)$, where
$G$ is finitely presented group (cf. the discussion in
~\cite{MR1085144}).

\subsection*{Properties of the fundamental group}
In Section \ref{S:properties_group}, we prove that the fundamental
group of a symplectically hyperbolic manifold is non-amenable.
In particular, it has exponential growth. We give a direct proof
of the last statement since it is essentially easier than the
proof of non-amenability which relies on a nontrivial result
of Gromov. Both results, however, follows from 
$\widetilde d$-boundedness of the volume form rather than symplectic
hyperbolicity.

\subsection*{Groups of diffeomorphisms}
In Section \ref{S:diff}, we prove that the group 
$\Ham\Mo$ of Hamiltonian diffeomorphisms of symplectically 
hyperbolic $\Mo$ is homotopy equivalent to $\Symp_0\Mo$,
the component of the identity of the group of all
symplectomorphisms. 

\subsection*{Further examples} 
\begin{enumerate}
\item
Let $F$ be a closed surface of genus at least two.
An oriented surface bundle $F\to E\to B$ over a symplectically
hyberbolic manifold admits a hyperbolic symplectic form
(see Corollary~\ref{C:surface}).
\item
If $F\to E\to B$ is a flat symplectic bundle, where the base
and the fibre are symplectically hyperbolic then the
total space admits a hyperbolic symplectic form
(see Therorem~\ref{T:flat_bundles}).
\item
Certain Lefschetz fibrations admits hyperbolic symplectic forms.
We give a precise statement and construction in Section 
\ref{SS:lefschetz}.
\end{enumerate}

\subsection*{Acknowledgements} 
I would like thank \'Swiatos\l aw Gal,
Dusa McDuff, Leonid Polterovich, 
Yuli Rudyak and Lukas Vokrinek 
for discussions and comments.

\section{Bounded classes are $\widetilde d$-bounded}
\label{S:bounded}

\begin{theorem}\label{T:bounded}
Let $M$ be a closed manifold. Let $\om\in \Om^k(M)$ be
a closed differential form such that its cohomology class
$[\om]\in H^k(M;\B R)$ is represented by a bounded cochain 
$c\in C^k(M,\B R)$. Then  $\om$ is $\widetilde d$-bounded.
\end{theorem}

\begin{proof}
Let $p:\widetilde M \to M$ be the universal cover.
We shall show that $p^*(\om)=d\alpha$,
where $\al$ is a form bounded with respect to
the metric $\widetilde g$ induced from a metric $g$ 
on $M$.

Let $c\in C^k(M,\B R)$ be a bounded singular cocycle representing
$[\om]$. Since $\widetilde M$ is simply connected, 
$p^*(c)= \delta (b)$ for some bounded real cochain $b\in C^{k-1}(M,\B R)$. 
Here we use the fact that if the fundamental group of
a space $X$ is amenable (e.g. solvable) then 
the bounded cohomology of $X$ is trivial in positive
degrees \cite{MR686042}.
Let us denote the bounding constant by $C\in \B R$.

Let $K$ be a smooth triangulation of $M$ and $K'$ the induced
triangulation of the universal cover $\widetilde M$. Let
$c'\in C^k(K,\B R)$ and $b'\in C^{k-1}(K',\B R)$  
be simplicial cochains induced 
by $c$ and $b$ respectively. We have that 
$\delta(b') = p^*(c')$ and $b'$ is bounded.
That is $\left <b',S\right> < C$ for any simplex  $S\in K'$.
Let $\psi_S$ be a cochain which attains the value $\pm 1$ 
on $\pm S$ and zero on other simplices.
We can express the cochain $b'$ as a sum
$b'=\sum_{S}r_S\psi_S$, where $S$ ranges over $(k-1)$-simplices
of $K'$ and $|r_S|<C$.

Let $\Phi_{K'}:C^*(K')\to \Om^*(\widetilde M)$ be a chain map which is a
right inverse to the integration over the simplicial chains
(see pages 148-149 in Singer-Thorpe \cite{MR0213982} for the 
definition and details). It also has the following property
(see Lemma 1 (4) on page 148 in \cite{MR0213982}):

\ms
\NI
{\em
If $S$ is an oriented simplex of the triangulation $K'$ 
then $\Phi_{K'}(\psi_S)$ is supported in the star of $S$.}

\ms
\NI
The construction of $\Phi_{K'}$ depends on the choice of
a partition of unity of $\widetilde M$ associated with the
triangulation $K'$. We take this partition of unity to be
induced from the partition of unity of $M$ used to define
$\Phi_{K}:C^*(K)\to \Omega^*(M)$. 
Hence it follows (from the construction of
~$\Phi_{K'}$) that if the triangulation $K$ of $M$ is fine 
enough then the maps $\Phi_{K'}$ commute with the
deck transormations. That is 
$\Phi_{K'}(h^*(\psi_S))=h^*(\Phi_{K'}(\psi_S)),$ where
$h:\widetilde M\to \widetilde M$ is a deck transormation.
In particular, we get the following 
$$
p^*\circ \Phi_{K} = \Phi_{K'} \circ p^*.
$$

Next we show that the form $\Phi_{K'}(b')=\sum_S r_S\Phi_{K'}(\psi_S)$ 
is bounded with respect to $\widetilde g$. First notice that
the  forms $\Phi_{K'}(\psi_S)$ are uniformly bounded. To see
this fix a fundamental domain in $\widetilde M$ for the action of 
$\pi_1(M)$. For any $S$ there exists $h\in \pi_1(M)$ such that
$h^*\psi_S= \psi_F$, where $F$ is a simplex in the fundamental domain.
Now the uniform boundedness follows from the fact that $\pi_1(M)$
acts on $(\widetilde M,\tilde g)$ by isometries.
Then, since $|r_S|< C$, we get that the differential form
$\Phi_{K'}(b')=\sum_S r_S\Phi(\psi_S)$ is bounded with respect
to $\widetilde g$ as claimed. 

We also claim that 
$d(\Phi_{K'}(b')) = p^*(\om + d\beta)$, for some form $\beta\in \Om^{k-1}(M).$
It is the following calculation.
\begin{eqnarray*}
d(\Phi_{K'}(b'))&=&\Phi_{K'}(\delta( b'))
=\Phi_{K'}(p^*(c'))=p^*(\Phi_{K}( c'))=
p^*(\omega + d\beta).
\end{eqnarray*}

\NI
Finally, we have that $p^*(\om)=d(\Phi_{K'}(b') - p^*(\beta))$ and 
$\al:=\Phi_{K'}(b') - d(p^*(\beta))$ is clearly bounded.

\end{proof}

\begin{corollary}\label{C:surface}
Let $F$ be a closed oriented surface of genus at least 2 and
let $(B,\om_B)$ be a symplectically hyperbolic manifold.
Then an oriented bundle $F\stackrel i\to M\stackrel p\to B$ 
admits a hyperbolic symplectic form.
\end{corollary}

\begin{proof}
Let $\om_F$ be an area form on $F$. 
The Thurston construction (Theorem 6.3 in \cite{MR2000g:53098})
gives a symplectic form in the class $C\cdot p^*[\om_B]+ \Om$, where
$\Om$ is any class in $H^2(M)$ such that $i^*(\Om)=[\om_F]$ and
$C>0$ is a constant large enough.

Let $\Om$ be a constant multiple of the Euler class of the
bundle $V:=\ker dp \to M$ tangent to the fibers of $p$.
According to Morita \cite{MR89c:57030} this class is bounded and so is
the class $p^*[\om_B]+\Om$.
\end{proof}

\section{Symplectic bundles and Lefschetz fibrations}\label{S:lefschetz}
In this section we give simple constructions of hyperbolic symplectic
forms on certain bundles and Lefschetz fibrations.

\subsection{Flat symplectic bundles}\label{SS:flat_bundles}
Recall that a symplectic form on the total space of a symplectic
bundle is called compatible if its pull-back to eaach fibre
is the symplectic form.

\begin{theorem}\label{T:flat_bundles}
Let $\Mo\to E\to B$ be a symplectic flat bundle. If
$\om$ and $\om_B$ are hyperbolic then $E$ admits a hyperbolic
compatible symplectic form.
\end{theorem}

\begin{proof}
Since the bundle is flat the total space is of the form
$E=(\widetilde B \times M)/\pi_1(B)$, where
the fundamental group of the base acts diagonally.
The form $\pi_{\tilde B}^*(\widetilde \om_B) + \pi_M^*(\om)$ 
is invariant under this action and it descend to a symplectic
from $\om_E $ on $E$. Hence on the universal cover $\widetilde E$
the induced form 
$\widetilde \om_E= \pi_{\widetilde B}^*(\widetilde \om_B)
+\pi_{\widetilde M}^*(\widetilde \om)$ is clearly a differential
of a bounded one-form.
\end{proof}

\subsection{Lefschetz fibrations}\label{SS:lefschetz}
In this section we construct a hyperbolic symplectic
form on certain Lefshetz fibrations. Let 
$p:E\to B$ be a 4-dimensional Lefschetz fibration.
According to Gompf and Thurston \cite[Theorem 10.2.18]{MR2000h:57038}, 
if there exists $\Om\in H^2(E)$ such that it restricts to
a non-trivial class of the fibre then the class
$\Om+C\cdot p^*[\om_B]$ admits a symplectic representative
for $C\in \B R$ big enough. Here $\om_B$ is the symplectic
form on the base. Thus if $\Om $ is $\tilde d$-bounded
and the base is of genus at least two then the symplectic
form is hyperbolic. In the next theorem we show that
we can construct such forms with some control
of the fundamental group. Let $\Pi_g$ denote the
fundamental group of a closed surface of genus $g$.

\begin{theorem}\label{T:lefschetz}
Let $\Om\in H^2(\Ga;\B R)$ be a bounded
class. 
For any $g\geq 2$
there exists a 4-dimensional closed symplectic manifold
$\Mo$ with $\pi_1(M)=\Ga\oplus \Pi_h$ such that the
symplectic form is hyperbolic.
\end{theorem}

\begin{proof}
The construction is the same as in K\c edra-Rudyak-Tralle 
\cite{SAS} and relies on the construction of
Amoros et al \cite{MR2002g:57051}.
First, we construct a symplectic Lefschetz fibration
$F\to X\to S^2$ such that:
\begin{enumerate}
\item
$\pi_1(X)=\Ga$.
\item
The pull-back $c^*(\Om)$ restricts nontrivially to the
fibre $F$; here $c:X\to K(\Ga,1)$ is the classifying map.
\item 
The fibration has a section.
\end{enumerate}
A construction of such Lefschetz fibration 
is provided in \cite{MR2002g:57051}.

Next take $F\times \Si_g$ equipped with the product symplectic structure.
Let $M:=X \#_F (F\times \Si_g)$ be the symplectic fibre sum.
Let $r:M \to X$ be a retraction. The cohomology class
$r^*(c^*(\Om))$ restricts nontrivially to the fibers hence the
Gompf-Thurston construction gives a symplectic form $\om$ in the
class $r^*(c^*(\Om)) + C\cdot p^*[\om_g]$. This class is
bounded hence the symplectic form $\om$ is hyperbolic.
The calculation of the fundamental group is a direct
application of the van-Kampen theorem.
\end{proof}

\begin{corollary}\label{C:lefschetz}
If $\Ga$ is a hyperbolic group with nontrivial $H^2(\Ga;\B R)$
then for any $g\geq 2$ there exists a closed symplectically
hyperbolic manifold $\Mo$ with $\pi_1(M)=\Ga\oplus \Pi_g$.\qed
\end{corollary}

\section{Isoperimetric inequality for symplectic forms}
\label{S:isoperimtric}

Let $M$ be a  manifold, $g$ a 
Riemannian metric and $\om\in \Om^2(M)$ a closed differential 
two-form.
We say that the form $\om$ satisfies the
{\bf linear isoperimetric inequality} (with respect to $g$) 
if there exists a constant
$C > 0$ such that for every smooth map $f:D^2\to M$ we have
$$
\left |\int_{f(D^2)}\om \right|\leq 
C\cdot \text{Length}(f(\partial D^2)).
$$
\bs

\begin{theorem}\label{T:isoperimetric}
Let $M$ be a closed  manifold.
A closed $2$-form $\om $ satisfies the linear isoperimetric inequality
if and only if it is $\widetilde d$-bounded.
\end{theorem}

\begin{proof}
Suppose that $\om$ is $\widetilde d$-bounded. Let $f:D^2\to M$ be a smooth
map. It admits a lift $\widetilde f:D^2\to \widetilde M$ to the
universal cover $p:\widetilde M\to M$. We calculate:
\begin{eqnarray*}
\left |\int_{D^2}f^*\om \right|&=&
\left |\int_{D^2}\widetilde f^*(p^*\om) \right|
= \left |\int_{D^2}\widetilde f^*d\alpha \right|\\
&=& \left |\int_{S^{1}}\widetilde f^*\alpha \right|
\leq C \cdot \text{Length}(\widetilde f(\partial D^2))\\
&=& C\cdot \text{Length}(f(\partial D^2)).
\end{eqnarray*}
Here, the constant $C$ bounds from above the norm of the
promitive $\alpha$, that is $\|\alpha(x)\|\leq C$ for
any $x\in \widetilde M$. Thus $\om$ satisfies an isoperimetric
inequality.

The converse is a direct application of the following result
due to Sikorav (Theorem 1.1 in \cite{MR1871291}).

\begin{theorem}[Sikorav]
Let $M$ be a Riemannian manifold with a triangulation of bounded
geometry. Let $\om \in \Om^q(M)$ be a closed form, and let
$f\in C^0(M,\B R_+)$ be such that
$$
\left |\int_T \om \right | \leq \int_{|\partial T|}f 
$$
for every simplicial chain $T\in C_q(K)$. Then $\om$
has a primitive $\alpha$ such that
$$
\al(x)\leq C_1 \max_{B(x,C_2)}(|\om| + f),
$$
for some constants $C_1,C_2$.\qed
\end{theorem}

\begin{remark}
We explain the notions in the formulation of the Sikorav
theorem (see \cite{MR1871291} for more details).
\begin{enumerate}
\item
The {\em triangulation of bounded geometry} is
a triangulation which satisfies the following two
properties:
\begin{enumerate}
\item
There exists a number $s\in \B N$ such that the link of
every simplex contains at most $s$ simplices.
\item
There exists a number $l>0$ such that for any
$q$-simplex $S$ of the triangulation there exist
a diffeomorphism $\psi_S:S\to \Delta^q$ such that
the norm of the differential of $\psi^{\pm 1}$ is bounded
by~$l$, $|d\psi_S^{\pm 1}|<l$. Moreover, $\psi$ can be extended (with
the norm condition preserved) to a neighborhood of $S$
sending it to a fixed neighborhood of 
$\De^{q}\subset \B R^n$, $n=\dim M$.
\end{enumerate}
Clearly, the triangulation of the unieversal cover of a closed
manifold has bounded geometry.
\item
If $S=\sum a_i \sigma_i$ is a $q$-chain, then
$$\int_{|S|}f := 
\sum_i |a_i|\int_{\Delta^q}\sigma_i^*(f\cdot \vol_q),$$
where $\vol_q$ is the $q$-dimensional volume induced
by the Riemannian metric.
\end{enumerate}
\end{remark}

\noindent
Now we can get back to the proof of Theorem \ref{T:isoperimetric}.
First observe that the isoperimetric inequality
immediately implies the asphericity of $\om$.
Moreover, the universal cover $(\widetilde M,\widetilde \om)$
also satisfies the isoperimetric inequality with the
same constant. We shall show that this 
implies the hypothesis of the Sikorav theorem. 

Let $K'$ be a smooth triangulation of $M$. Then the triangulation
$K$ induced on the universal cover $\widetilde M$ has
bounded geometry. Let $T\in C_2(K)$ be a simplicial chain.
Its boundary is a cycle so it is a sum of loops
$\partial T = \sum_i\gamma_i$. Accroding to the
simple connectivity of $\widetilde M$ there exists
a chain $T_i$ such that $\partial T_i=\gamma_i$ and
$T_i$ is the image of a triangulation of the disc $D^2$.
Since the sum $T - \sum_iT_i$ is a cycle its symplectic
area is zero, $\int_{T - \sum T_i}\tilde\omega$. This is true
because $\widetilde M$ is simply connected so any homology
class of degree two is represented by a sum of
spheres and $\tilde \omega $ is aspherical.

In other words, we have
$\int_{\sum T_i} \widetilde \om = \int _T\widetilde \om$
and we calculate
\begin{eqnarray*}
\left |\int_T\widetilde \om\,\right| &=&
\left |\int_{\sum T_i}\widetilde \om\,\right| 
\leq
\sum \left |\int_{T_i}\widetilde \om\,\right| \leq \\
&\leq& 
\sum C\cdot \text{Length}(\gamma_i)
=
\int _{|\partial \sum T_i|}C = \int _{|\partial T|}C 
\end{eqnarray*}

\NI
Now according to Sikorav's theorem there exists a primitive
$\al$ such that
$$
|\al(x)|\leq C_1 \max_{B(x,C_2)}(\,|\om| + C\,).
$$
That is $\al$ is bounded since $\om$ is bounded.

\end{proof}

\section{Symplectic hyperbolicity and the fundamental group}
\label{S:group}
In this section we prove that the symplectic hyperbolicity
depends in a sense on the fundamental group only. More precisely,
we have the following result.

\begin{theorem}\label{T:group}
Let $\Mo$ be a closed symplectically hyperbolic manifold with 
$[\om] = f_M^*(\Om)$, where $\Om\in H^2(\pi_1(M);\B R)$.
If $\Wo$ is a symplectic manifold with $\pi_1(W)=\pi_1(M)$
and $[\om_W]=f_W^*(\Om)$ then $\Wo$ is symplectically hyperbolic.
\end{theorem}

\begin{proof}
Let $K_n$ denotes the $n$-skeleton of $K(\pi_1(M),1)$. 

\begin{lemma}\label{L:main}
There exists a manifold $M_n$ such that 
$K_n\subset M_n\stackrel{i}\hookrightarrow K(\pi_1(M),1)$
and $i^*(\Om)$ is hyperbolic.
\end{lemma}

The classifying map $f_W:W\to K(\pi_1(W),1)$ factorizes through
some finite dimensional skeleton $K_n$ hence through $M_n$ as well.
The hyperbolicity of $[\om_W]$ follows since $i^*(\Om)$ is hyperbolic.

\end{proof}

\begin{proof}[Proof of Lemma \ref{L:main}]
Regard $K(\pi_1(M),1)$ constructed from $M$ by attaching cells.
We want to perform this construction so that at every stage it
is a manifold. This is done as follows. Let $\si:S^k\to M$ represent
a homotopy class $[\si]\in \pi_k(M)$ we are going to kill. Take a product
$M\times D^m$ so that $[\si]$ is represented by an embedding
$s:S^k\to M\times S^{m-1}=\partial (M\times D^m)$.
Then attach a handle $D^{k+1}\times D^{2n+m-k}$ along $s$ (after choosing
any framing of $s$). Call the resulting manifold $M_{\si}$.
The classifying map $f_M:M\to K(\pi_1(M),1)$ clearly extends
to $M_{\si}$. We shall show that the class $f_{M_{\si}}^*(\Om)$
is hyperbolic.

Denote by $\om\in \Om^2(M\times D^m)$ the pull-back of the symplectic 
form on $M$ under the projection. This form extends to a form $\om_{\si}$
on $M_{\si}$ representing $f_{M_{\si}}^*(\Om)$. 
Choose a Riemannian metric $g$ on 
$M\times D^m$ so that $\widetilde \om = d\al$ on the universal 
covering $\widetilde M\times D^m$ and $\al$ is bounded with respect
the the metric induced from $g$. 
Extend the metric $g$ to $g_{\si}$ on $M_{\si}$.
We need to show that the induced form
$\widetilde \om_{\si}$ on the universal cover of $M_{\si}$
is a differential of a bounded one-form.

We have that $\widetilde \om_{\si}=d\be$ for some one-form $\beta$.
The universal cover $\widetilde M_{\si}$ is $\widetilde M\times D^m$ with
infinitely handles $H_{\ga}:=D^k\times D^{2n+m-k}$ attached. The handles correspond
to elements $\ga$ of the fundamental group $\pi_1(M)$. The induced form 
$\widetilde \om_{\si}$ is the same when restricted to every handle.
Choose its primitive $\widehat \al$. It is bounded due to the compactness
of the handle.

Let $(R,S)$, where $R\in [0,1]$ and $S\in S^k$
be polar coordinates on $D^{k+1}$. For $R\in [1-\varepsilon,1]$
we have $\widetilde \om_{\si}=d\al=d\widehat \al$ on the handle
$H_{\ga}$. In this neighborhood we take the covex combination
(in the affine space $\al + \ker d$) interpolating between them. 
This operation preserves the boundedness. Hence we obtain
a bounded primitive of $\widetilde \om$.

The rest of the proof is to apply the induction on the
handles attached in order to obtain $K(\pi_1(M),1)$ from
$M$. This finishes the proof of Lemma \ref{L:main}.
\end{proof}

Recall that a $\pi_1$-cobordism between manifolds $M_1$ and
$M_2$ is a cobordism $W$ such that the inclusions induce
isomorphisms on the fundamental group.

\begin{corollary}\label{C:cobordism}
Let $(M_1,\om_1)$ and $(M_2,\om_2)$ be closed symplectic
manifolds and $W$ a $\pi_1$-cobordism between them. Suppose that
there exist $\Om\in H^2(W;\B R)$ such that it pulls back to the
classes of the symplectic forms under the inclusions
$\iota_i:M_i\to W$. 
Then $\om_1$ is hyperbolic if and only
if $\om_2$ is hyperbolic.
\end{corollary}

\section{Properties of the fundamental group of a symplectically
hyperbolic manifold}\label{S:properties_group}

\subsection{Infiniteness}\label{SS:infty}
\begin{proposition}\label{P:infty}
Let $\Mo$ be a closed symplectically hyperbolic
manifold. Then
\begin{enumerate}
\item
The fundamental group of a $\Mo$ is infinite.
\item
There exists $\Om\in H^2(K(\pi_1(M));\B R)$ such that
$c_M^*(\Om)=[\om]$, where $c_M:M\to K(\pi_1(M),1)$
is the classifying map.
\item
If $h:\B Z\oplus \B Z\to \pi_1(M)$ is a homomorphism
then $h^*\Om=0$
\end{enumerate}
\end{proposition}

\begin{proof}\hfill
\begin{enumerate}
\item
If $\pi_1(M)$ were finite the the universal cover $\widetilde M$
would be closed and hence the induced symplectic form would
not vanish on some sphere. The same would hold for $\om$ which
contradicts Proposition \ref{P:properties}.
\item
Consider the homotopy fibration 
$$\widetilde M \stackrel{p}\to M\stackrel{c_M}\to K(\pi_1(M),1)$$
and the associated spectral sequence.
Since $\widetilde M$ is simply connected and $p^*[\om]=0$,
we get that $[\om]=c_M^*(\Om)$ for some 
$\Om\in H^2(K(\pi_1(M),1),\B R)$.
\item
Suppose that $h:\B Z^2\to \pi_1(M)$ is a homomorphism. It induces
a continuous map $h':T^2=K(\B Z^2,1)\to K(\pi_1(M),1)$. Since
$\widetilde M$ is simply connected and $\dim T^2=2$, the map
$h'$ admits a lift to $H:T^2\to M$. We have that
$c_M\circ H=h'$. Due to Proposition~\ref{P:properties}, we
know that the symplectic form vanishes on tori. Thus
$0=H^*[\om]=H^*(c_M^*\Om)= h^*(\Om)$ which finishes the proof.  
\end{enumerate}
\end{proof}

\begin{remark}
The first two statements hold for aspherical symplectic forms.
\end{remark}

\subsection{Exponential growth}\label{SS:growth}

\begin{proposition}
Let $\Mo$ be closed symplectic manifold. If the symplectic
form $\om$ is hyperbolic then the fundamental group $\pi_1(M)$
has exponential growth.
\end{proposition}

\begin{proof}
It is known that the fundamental group of a manifold has
exponential growth if and only if balls in the universal
cover grow exponentially with respect to the radius.
We shall prove the latter.

Let $(\widetilde M,\widetilde \omega)$ be the universal cover
of $\Mo$ and let $\widetilde g$ be the Riemannian metric induced
from a metric $g$ compatible with $\om$. Let $\al\in \Om^1(\widetilde M)$
be a primitive of $\widetilde \om$ and let $X$ be the corresponding
vector field. That is $\iota_X\widetilde \om = \al$. Since
$\widetilde g$ is compatible with $\widetilde \om$, the norm
of $X$ is uniformly bounded by the constant bounding the
norm of $\al$. Let $\vol=\widetilde \om^n$ denote the volume
form. We have that
$
L_X\vol = L_X\widetilde \om^n = n\cdot \widetilde \om^n =
n\cdot \vol.
$

Let $B:=B(x,1)$ denote the ball of radius 1 with center at $x\in \widetilde M$.
Let $\psi:\B R\to \Diff(M)$ be the flow corresponding to the
vector field $X$. We claim that the volume of the image
$\psi_t(B)$ grows exponentially with $t$. This is the
following calculation.
\begin{eqnarray*}
\frac{d}{dt}\Big|_{t=s}\vol(\psi_t(B))&=& 
\frac{d}{dt}\Big|_{t=s}\int_{\psi_t(B)}\vol\\
&=&
\frac{d}{dt}\Big|_{t=s}\int_{B}\psi_t^*\vol\\
&=&
\int_{B}\frac{d}{dt}\Big|_{t=s}\psi_t^*\vol\\
&=&
\int_{B}\psi_s^*(L_X\vol)\\
&=&
\int_{\psi_s(B)}n\cdot\vol=n\cdot \vol(\psi_s(B))
\end{eqnarray*}
Hence we get that $\vol(\psi_t(B))=e^{nt}\cdot \vol(B)$.

On the other hand $\psi_t(B)\subset B(\psi_t(x),2\cdot C\cdot t + 1)$ 
which implies that 
$\vol\left(B\left(\psi_t(x),2Ct+1\right )\right )\geq e^{nt}\vol(B(x,1))$.
Since the deck trasformations are isometries of $\widetilde g$,
we get that
$$\vol\left(B\left(\psi_t(x),2Ct+\text{diam}(M)+1\right)\right)
\geq e^{nt}\vol(B(x,1)).$$
The diameter constant is added because the deck transormation
might not move $\psi_t(x)$ to $x$. This proves that the volume
of balls centered at $x$ grow exponentially with respect to the
radius and hence it proves the exponential growth of $\pi_1(M)$.
\end{proof}

\subsection{Non-amenability}\label{SS:amenable}

\begin{theorem}\label{T:amenable}
Let $\Mo$ be a closed symplectically hyperbolic manifold.
Then its fundamental group $\pi_1(M)$ is non-amenable.
\end{theorem}

\begin{proof}
Let $\al\in \Om^1(\widetilde M)$ be a bounded primitive of
$\widetilde \om$. Then the form $\al\wedge \widetilde \om\,^{n-1}$
is a bounded primitive of the volume form $\widetilde \om\,^n$.
This implies that $(\widetilde M,\widetilde g)$ satisfies
an isomerimetric inequality:
$$
\int_{X}\widetilde\om\,^n=\int_{\partial X}\al\wedge \widetilde \om\,^{\,n-1}
\leq C\cdot \vol_{2n-1}(\partial X),
$$
where $X\subset \widetilde M$ is any domain
and $\vol_{2n-1}$ denotes the $2n-1$-dimensional volume induced
by the Riemannian metric $\widetilde g$.
It follows from Theorem 6.19 in \cite{MR2307192} 
that $\pi_1(M)$ is not amenable.
\end{proof}

\begin{example}\label{E:amenable}
Let $\Mo=G/\Gamma$ be a closed symplectic solvmanifold.
That is it is a homogeneous space of a simply connected
solvable Lie group $G$ and $\pi_1(M)=\Gamma\subset G$ 
is a lattice. Since $\Gamma$ is amenable, $\om$ is
not hyperbolic. If $\Gamma$ does not contain a nilpotent
subgroup of finite index then it as exponential growth.
\end{example}

\section{Groups of diffeomorphisms}\label{S:diff}

Recall that
the {\bf flux group} $\Gamma_{\om}$ is the image of 
the flux homomorphism $$\Flux:\pi_1(\Symp_0\Mo)\to H^1(M;\B R).$$ 
It is defined
by $$\left<\Flux[\xi_t],[A]\right>= \int _{T^2}\xi_A^*\om,$$
where $\xi_A(s,t)= (\xi_t)(A(s))$ and 
$\xi:S^1\to \Symp_0\Mo$ is a loop based at the identity and
$A:S^1\to M$ is a 1-cycle in $M$ 
(see Section 10.2 in \cite{MR2000g:53098} for more details).
Notice that if the flux homomorphism is non-trivial then
the symplectic for does not vanish on some torus.
The groups of symplectic and Hamiltonian diffeomorphisms
form the following extension (Theorem 10.18 in \cite{MR2000g:53098}).
$$\Ham\Mo\to \Symp_0\Mo \to H^1(M;\B R)/\Gamma_{\om}$$
The next proposition is an immediate consequence of
the definition of the flux homomorphism.

\begin{proposition}\label{P:diff}
Let $\Mo$ be a closed symplectic manifold such that
the symplectic form vanishes on tori.
Then $\Symp_0\Mo\simeq\Ham\Mo$. In other words, the flux
group of $\Mo$ is trivial. In particular, the statement
holds for symplectically hyperbolic $\Mo$, due to
Proposition \ref{P:properties}.
\qed
\end{proposition}

Polterovich proved (among others very interesting and
beautiful results) in \cite[1.6.C]{MR2003i:53126} that if
$G\subset \Ham\Mo$ is a finitely generated subgroup
and $\Mo$ is symplectically aspherical then cyclic
subgroups are undistorted in $G$ with respect to the
word metric. More precisely, if $\text{Id}\neq g\in G$
then $\|g^n\|\geq C\cdot |n|$ for some $C>0$,
where $\|g\|$ is the norm (the distance from the identity) 
given by the fixed finite set of generators. The property of
being undistorted does not depend on the choice 
of a finite set of generators.

\begin{example}\label{E:dist}\hfill
\begin{enumerate}
\item
Any cyclic subgroup of a free or free abelian group is undistorted.
Indeed, it is easy to check that $\|g^n\|=n\cdot\|g\|$ in this case.
\item
Torsion groups are not undistorted. 
\item
Let $1<p<q$ be integers. The subgroup of the Baumslag-Solitar group
$G:=\left <x,t\,|\,x^q=tx^pt^{-1}\,\right>$ generated by
$x$ is not undistorted.
\item
Let  $h:G\to \B Z^k$ be a surjective homomorphism. If
$h(g)~\neq~0$ then $g$ generates an undistorted subgroup in $G$.
Indeed, let $g_1,\dots,g_n\in G$ be generators of $G$ such that
$h(g_1),\dots,h(g_k)$ are the standard generators of $\B Z^k$. 
With respect to this sets of generators, we have that 
$\|h(g)\|\leq C\cdot \|g\|$ for some $C>0$. 
Since $\|(h(g))^n\|=|n|\cdot\|h(g)\|$, we get that 
$\|g^n\|\geq c\cdot |n|$, for some $c>0$.
\end{enumerate}
\end{example}

The next proposition is a slight generalization
of the above mentioned Polterovich's result.

\begin{proposition}\label{P:polt}
Let $\Mo$ be a symplectically hyperbolic manifold.
If $G\subset \Symp_0\Mo$ is a finitely generated subgroup
then its cyclic subgroups are undistorted.
\end{proposition}

\begin{proof}
If $G\subset \Symp_0\Mo$ is a subgroup then
we have a morphism of extensions (the vertical
arrows are inclusions):
$$
\xymatrix
{
K \ar[r]\ar[d] &  G\ar[r]\ar[d] &  Q\ar[d]\\ 
\Ham\Mo\ar[r] & \Symp_0\Mo \ar[r] & H^1(M;\B R)
}
$$
Let $\text{Id}\neq g\in G$. If it belongs to $K$ then
according to Polterovich's theorem it generates
an undistorted subgroup. If it maps non-trivially
to $Q$ then it also generates an undistorted subgroup
because $Q$ is free abelian (see Example \ref{E:dist} (4)).
\end{proof}

\bibliography{../../bib/bibliography.bib}
\bibliographystyle{plain}

\end{document}